%% file: arxiv_submission.tex
\documentclass[]{amsart}

\input{preamble.tex}

\begin{document}

\title[Steady generalized Ricci solitons]{
  Steady Generalized Ricci Solitons with
  $\mathrm{SO}(3)$ Symmetries in Dimension Three
}

\author[E. Cochran]{Eric Cochran}
\author[M. Womack]{Michael Womack}

\begin{abstract}
    We completely classify all gradient, complete, steady generalized Ricci solitons in dimension three with $SO(3)$ symmetry, extending a result of Podest\'a-Raffero. We further extend this analysis to $n$-forms on $n$-manifolds with an $SO(n)$ symmetry.
\end{abstract}

\maketitle

\section{Introduction}

Let $(M,g)$ be a Riemannian manifold, and let $H$ be a closed 3-form on $M$.  The pair $(g,H)$ is said to be a \textit{steady generalized Ricci soliton} if there exists a vector field $X \in \mathfrak{X}(M)$ which satisfies 
\begin{equation}\label{eq:1}
\begin{split}
    \Ric - \frac{1}{4}H^2 + \frac{1}{2}\mathcal{L}_X g = 0, \\
    \frac{1}{2}d^*H  + \frac{1}{2}i_XH  = 0.
\end{split}
\end{equation}
The soliton is said to be \textit{gradient} if $X=\nabla f$ for some $f\in C^{\infty}(M)$. In a similar vein to how Ricci solitons are self-similar solutions of the Ricci flow, generalized Ricci solitons are self-similar solutions of the \textit{generalized Ricci flow}.  This flow is defined by the coupled equations
\begin{align} \label{eq:flow}
    &\frac{\partial }{\partial t} g = -2\mathrm{Ric}+\frac{1}{2}H^2
    \\&\frac{\partial }{\partial t} H = -\Delta_{d,g_t}H
\end{align}

\noindent where $H^2(X,Y) := 2\langle i_XH, i_YH \rangle$ and $\Delta_{d,g_t}$ denotes the Hodge Laplacian.  

In this paper, we completely classify steady, gradient generalized Ricci solitons in dimension three with $SO(3)$-symmetry. Specifically, we prove the following:

\begin{thm} \label{thm:main}
    For every $\lambda \in \mathbb{R}$, there exists a steady, rotationally symmetric, generalized Ricci soliton $(g_\lambda, H_\lambda)$ in dimension three. Moreover,
        \begin{itemize}
            \item[(a)] If $\lambda < 0$, the soliton $(g_\lambda, H_\lambda)$ is complete and noncompact.
            \item[(b)] If $\lambda = 0$, the soliton $(g_\lambda, H_\lambda)$ is complete and asymptotically conical.
            \item[(c)] If $\lambda > 0$ and $\lambda \neq 2$, the metric encounters a singularity at finite distance and is therefore incomplete.
             \item[(d)] If $\lambda = 2$, the corresponding soliton $(g_{\lambda}, H_\lambda)$ is complete, and is the Bismut flat 3-sphere.
        \end{itemize}
\end{thm}

\begin{thm}
    Up to homothety, every complete, steady generalized Ricci soliton with $H \not\equiv 0$ admitting an effective, nontransitive $SO(3)$ action is a member of the classification of Theorem \ref{thm:main} or a quotient of the Bismut flat $S^3$.
\end{thm}

The parameter $\lambda$ in Theorem \ref{thm:main} has a geometric interpretation. It is called the \textit{dilaton constant}, and is a conserved quantity of generalized Ricci solitons. See Proposition 1.1 of \cite{Podesta25}.  The parameter $\lambda$ in this paper is related to the parameter $q$ in the Podest\`a-Raffero paper by $\lambda = 6q + 5$. An analogous parameter and conserved quantity was used to study and classify steady gradient Ricci-Yang-Mills solitons on surfaces in \cite{Womack26}.

Steady gradient Ricci solitons in the classical case in dimension three are well understood. In \cite{Bryant}, Bryant constructs a one parameter family of Ricci solitons in dimension three with $SO(3)$-symmetries, and furthermore shows there is a unique steady example which is complete and has positive curvature. In later work, Brendle showed that a steady, gradient, complete, non-flat Ricci soliton in dimension three which is asymptotic to the Bryant soliton must be isometric to the Bryant soliton (see \cite{Brendle13}).

In recent years, steady generalized gradient Ricci solitons in dimension three have also been studied. Carmona studies such solitons with $SO(2) \times \mathbb{R}$ symmetry in \cite{Carmona26}. In \cite{Podesta25}, Podest\`a and Raffero construct complete Bryant-like solitons in the generalized setting. Specifically, they find a one parameter family of noncompact, steady generalized Ricci solitons in dimension three with positive curvature. Theorem \ref{thm:main} extends this result, as we identify all $\lambda$ values for which the associated soliton is complete.

We further consider analogous results for $SO(n)$-invariant solutions of a coupling of the Ricci flow equations to the evolution of a closed $n$-form. In this setting, many of the results from the dimension three case extend easily, so we are able to construct higher dimensional generalized Bryant-type steady solitons. The corresponding classification result does not extend as cleanly in higher dimensions however, for reasons described at the end of section two.

\subsection*{Acknowledgments:} Both authors would like to thank their research supervisor, Jeff Streets, for his helpful insights, guidance, and comments. Both authors were partially supported by NSF grant DMS-2342135.

\subsection*{AI Use Disclosure:} The main proof techniques in this paper were discovered with the assistance of ChatGPT 6-Astra and Claude Opus 5.5. The exposition was written entirely by the authors, and the proofs were rewritten by the authors for improved readability. The authors checked the proofs and take full responsibility for the correctness of the results.

\section{Results}

We split up the proof of Theorem 1.1 into several smaller results. Before tackling the proof, we begin by defining some notation.

First, following \cite{Podesta25}, we write the metric as
\begin{align}
    g = dt^2 + \phi^2(t)g_{S^2}
\end{align}
\noindent which is the general metric for a three dimensional metric with $SO(3)$ symmetry.

We now define the parameter $\lambda$ as follows:

\begin{align}
    \lambda := f''+2\frac{\phi'}{\phi}f'+2e^{2f}-(f')^2
\end{align}

\noindent which, as mentioned earlier, is a conserved quantity shown in \cite{Podesta25}.

Next, since we are working in dimension three, the closed three form $H$ is determined by a single function.  Considering a generic $SO(3)$-invariant closed three form $H$ and gradient vector field $X = \nabla f$ and plugging in $H$, $g$, and $X$ into \ref{eq:1} (following \cite{Podesta25}), we obtain $H(t) = k\phi^2e^f dt \wedge e^{12}$ for some constant $k$ as determined in equation (3.8) of \cite{Podesta25}. Throughout this paper, we make the normalization $k=\sqrt2$. With this normalization, we further obtain the following system of ODEs:
\begin{equation} \label{eq:phi_f_sys}
\begin{split}
    \phi''  &=\frac{1-(\phi')^2}{\phi}+\phi'f'-\phi e^{2f}, \\
    f'' &= 2\frac{1-(\phi')^2}{\phi^2}+2\frac{\phi'}{\phi}f'-e^{2f}.
\end{split}
\end{equation}

Initial conditions for $\phi$ and $f$ are determined by the requirement that the solution extends smoothly across the origin, hence
\[
    \phi(0)=0, \qquad \phi'(0)=1,\qquad f(0)=f'(0)=0.
\]
Note that $f$ is only determined up to adding a constant, so we choose $f(0)=0$. Near the origin, the solutions $\phi$ and $f$ satisfy
\begin{equation}\label{eq:orig_asymp}
    \phi(t) = t+\frac{\lambda-5}{36}t^3+O(t^5),\qquad f(t)=\frac{\lambda-2}{6}t^2+O(t^4).
\end{equation}

\noindent as shown by Podest\`a-Raffero in page of 11 of their paper \cite{Podesta25}.

\begin{proposition}
    Suppose $(M,g,H,X)$ satisfies equation \eqref{eq:1} in dimension three with $H\neq 0$. If $M$, $g$ and $H$ are invariant under an $SO(3)$ action, there exists an $SO(3)$-invariant smooth function $f$ satisfying
    \[
        H=\sqrt{2}e^f\;dV_g,\qquad X=\nabla f.
    \]
\end{proposition}
\begin{proof}
    In dimension three we may write $H=q\; dV_g$ for some $q\in C^\infty(M)$. The second equation in \eqref{eq:1} is $d^*_g(q\;dV_g)=-i_{\nabla q}dV_g$, so
    \[
        dq=qX^\flat.
    \]
    From this, along every smooth curve $\sigma$ we see that the function $q\circ\sigma$ solves the linear ODE
    \[
        (q\circ\sigma)' = X^\flat(\dot\sigma)(q\circ\sigma)
    \]
    and thus by uniqueness for solutions to this equation, if $q$ vanishes at a point then it must vanish everywhere. We can now assume $q>0$ everywhere, and setting 
    \[
        f=\log \frac{q}{\sqrt2}
    \]
    we have $df=X^\flat$.
\end{proof}
\begin{rmk}
    Compact generalized Ricci solitons were already known to be gradient (see \cite{garcia-fernandezGeneralizedRicciFlow2021}). The above proposition shows that this automatically true in the noncompact setting assuming $SO(3)$-invariance.
\end{rmk}

\subsection{Subcritical/noncompact solutions} We begin by studying case (a) of Theorem \ref{thm:main}, where we assume $\lambda <0$, and prove the long-time existence of the system of ODEs (\ref{eq:phi_f_sys}). 

A key quantity which will prove to be useful throughout this paper is the weighted mean curvature of the measure space $(M, g, e^{-f} \mathrm{dV_g})$. For the ansatz we are making above, the weighted mean curvature is

\[ h  = 2 \frac{\phi'}{\phi} - f'. \]

We begin by proving a useful lemma, which we will use to establish the long-time existence of the solution to the system \eqref{eq:phi_f_sys} for $\lambda \leq 0$.

\begin{lem}
    The quantity $h$ is positive and decreasing in $t$ for $\lambda \leq 0$.
\end{lem}

\begin{proof}
    To show $h > 0$, we first compute

    \[ h^2 = \left ( 2 \frac{\phi'}{\phi} - f'\right )^2 = 4\frac{(\phi')^2}{\phi^2}-4f'\frac{\phi'}{\phi}+(f')^2. \]
    From here, we recall an important conservation law for generalized Ricci solitons.  From Proposition 1.1 of \cite{Podesta25}, we have 
        \[ f'' + 2\frac{\phi'}{\phi}f' + 2e^{2f} - (f')^2= C \]
    for some constant $C$.  Using the second equation of (\ref{eq:phi_f_sys}), we can rewrite this conservation law as 
    \begin{equation}\label{eq:cons_law}
        2 \frac{1-(\phi')^2}{\phi^2} + 4\frac{\phi'}{\phi}f' - (f')^2 +  e^{2f} = \lambda.
    \end{equation}
    Using this conservation law, we find that

    \begin{align} \label{eq:h_square}
     h^2 = 2 \left ( \frac{\phi'}{\phi} \right )^2+\frac{2}{\phi^2} + e^{2f} - \lambda. 
    \end{align}
    From this expression, we see that $h \neq 0$.  Furthermore, for small values of $t$, as $t \to 0$, we have that
    \[
        h(t) = \frac{2}{t} + O(t)  
    \]

    \noindent so that near zero, $h > 0$, and so $h(t)>0$ on its maximal interval of existence.

    To show that $h$ is decreasing, we compute
    \begin{equation}\label{eq:Hf_deriv}
        h' = 2\frac{\phi''}{\phi} - 2\left ( \frac{\phi'}{\phi}\right)^2 - f'' = -2 \left (\frac{\phi'}{\phi} \right )^2 - e^{2f} < 0
    \end{equation}
    \noindent where the second equality comes from using the system (\ref{eq:phi_f_sys}) to rewrite $h'$.  
\end{proof}

We now proceed with the following, showing global existence of the solution of the system (\ref{eq:phi_f_sys}) for $\lambda \leq 0$, implying part (a) of Theorem \ref{thm:main}.

\begin{prop}\label{prop:glob_exist}
    If $\lambda \leq 0$, then the local solution to the system of ODEs (\ref{eq:phi_f_sys}) extends to a global solution on $[0, \infty)$.
\end{prop}

\begin{proof}
    By contradiction, suppose that the solution does not exist for all time, and let $(0, L]$ be the maximal interval of existence where $0<L<\infty$.  Fix $t_0 \in (0, L)$.  By the previous lemma, $h(t)$ is positive and decreasing, and so
        \[ 0 < h(t) < h(t_0) \]
    for $t \in (t_0, L]$.  Using the expression for $h^2$ derived in the previous lemma, we have the following:
    \[
        \left |\frac{\phi'}{\phi} \right | \leq \frac{1}{\sqrt{2}} h(t_0), \qquad \frac{1}{\phi} \leq \frac{1}{\sqrt{2}} h(t_0), \qquad e^f \leq h(t_0).
    \]
    From this, we have uniform bounds for $\phi, \phi', f, f'$ and $\phi^{-1}$.  Equation \ref{eq:phi_f_sys} then also gives uniform bounds on $\phi''$ and $f''$.  By standard ODE existence theory, since the solution has a finite limit as $t$ approaches $L$, the local solution can be extended beyond $L$, and so the solution must exist for all time as desired.
\end{proof}

\subsection{Critical/conical solutions}
Now, we turn to the case when $\lambda = 0$, and show that the complete soliton is asymptotically conical.

\begin{lem}\label{lem:crit_h_phi_asymp}
    When $\lambda = 0$, as $t\to\infty$ we have
    \[
        h \to 0, \qquad \phi\to\infty.
    \]
\end{lem}

\begin{proof}
    We have previously shown that the solution exists for all $t\geq 0$ and that $h$ is positive and decreasing. This implies that $h$ decreases to some fixed $\kappa\geq 0$ as $t\to \infty$. 

    Suppose $\kappa>0$. First note that by \eqref{eq:h_square} and the fact that $h(t)$ is positive and decreasing, $\frac{\phi'}{\phi}$ and $e^{f}$ are uniformly continuous on $[t_0, \infty)$. Integrating \eqref{eq:Hf_deriv}, we see that
    \[
        \int_{t_0}^\infty \left(2\left(\frac{\phi'}{\phi}\right)^2+e^{2f}\right)dt <\infty.
    \]
    Hence $\frac{\phi'}{\phi}\to 0$ and $e^{f}\to 0$ as $t\to \infty$. The conservation law then implies
    \[
        \lim_{t\to\infty} h^2=\lim_{t\to\infty} \left(2\left(\frac{\phi'}{\phi}\right)^2+\frac{2}{\phi^2}+e^{2f}\right) =\lim_{t\to\infty} \frac{2}{\phi^2} = \kappa^2 > 0
    \]
    so $\frac{1}{\phi^{2}}\to \frac{\kappa^2}{2}$. However, this implies that
    \[
        \left(\frac{\phi'}{\phi}\right)' = \frac{1}{\phi^2}-2\left(\frac{\phi'}{\phi}\right)^2+\frac{\phi'}{\phi}f'-e^{2f}\to \frac{\kappa^2}{2}>0
    \]
    which contradicts that $\frac{\phi'}{\phi}\to 0$. Thus we must have $\kappa=0$, so $\phi\to\infty$ and $h\to 0$.
\end{proof}

\begin{lem}
    When $\lambda = 0$, we have that $f'<0$ on $(0,\infty)$.
\end{lem}

\begin{proof}
    From equation \eqref{eq:orig_asymp} we have that near the origin for the critical solution, 
    \[
        f'(t) = -\frac{2}{3}t+O(t^3)
    \]
    so $f'<0$ near the origin. We can combine equations \eqref{eq:phi_f_sys} and \eqref{eq:cons_law} to yield that
    \[
        f''=-h f'-2e^{2f}
    \]
    so that at a zero of $f'$, $f''<0$. This means $f'$ may only cross zero from $f'\geq 0$ to $f'<0$, and thus $f'$ remains below zero.
\end{proof}

We now derive more explicit asymptotics for $\phi$ and $f$. A key ingredient here is the conservation law \eqref{eq:cons_law}, which allows us to reduce the soliton system of equations to a coupled planar system of ODEs and study its trajectories (see \cite{perkoDifferentialEquationsDynamical2001}). We start by considering a set of substitutions to normalize the quantities in the conservation law. Let
\[
    a=\frac{\phi'}{\phi}, \qquad b=\frac{1}{\phi}, \qquad c=e^f
\]
so that the conservation law becomes
\[
    h^2 = 2a^2+2b^2+c^2
\]
and normalizing by $h$ defines
\[
    X=\frac{a}{h}, \qquad Y=\frac{b}{h}, \qquad Z=\frac{c}{h}, \qquad \frac{ds}{dt}=h.
\]
The conservation law becomes
\begin{equation}\label{eq:norm_cons_law}
    2X^2+2Y^2+Z^2=1,
\end{equation}
and the above quantities evolve by
\begin{align*}
    \dot{X} &= Y^2-Z^2-X+X(2X^2+Z^2),\\
    \dot{Y} &= Y(2X^2+Z^2-X),\\
    \dot{Z} &= Z(2X-1+2X^2+Z^2),
\end{align*}
where a dot denotes $d/ds$. Since $\phi\to\infty$ and $|a|\leq h/\sqrt{2}$,
\[
    \log\frac{\phi(t)}{\phi(t_0)} = \int_{t_0}^{t} a \;du\leq \frac{1}{\sqrt{2}}\int_{t_0}^{t} h\;du
\]
so $s\to\infty$ as $t\to\infty$.

Note that $Y\ne 0$, and by applying \eqref{eq:norm_cons_law} to eliminate $Y$ from the evolution of $X$, we reduce the system to a planar system of ODEs on the open half-disk
\begin{equation}\label{eq:planar_region}
    \mathcal{D} = \{(X,Z)\,:\, 2X^2+Z^2<1,\; Z>0\}
\end{equation}
with corresponding evolution equations
\begin{equation}\label{eq:norm_system}
\begin{split}
    \dot{X} &= \frac{1}{2}(4X^3-2X^2+2XZ^2-2X-3Z^2+1),\\
    \dot{Z} &= Z(2X^2+2X+Z^2-1).
\end{split}
\end{equation}

\begin{lem}\label{lem:crit_asymptotics}
    In the region $\mathcal{D}$, the normalized critical trajectory converges to
    \[
        (X,Y,Z)\to \left(\frac{1}{3},\frac{1}{\sqrt{3}},\frac{1}{3}\right).
    \]
\end{lem}

\begin{proof}
    By the Poinecar\'e-Bendixon Theorem, a trajectory of \eqref{eq:norm_system} must converge to either a critical point or a periodic orbit of \eqref{eq:norm_system} in $\overline{\mathcal{D}}$. The relevant critical points in the closure of $\mathcal{D}$ are
    \[
        \vec{x}_1 = \left(\frac{1}{2},\frac{1}{2},0\right),\quad \vec{x}_2 =\left(\frac{1}{\sqrt{2}},0,0\right),\quad \vec{x}_3 =\left(-\frac{1}{\sqrt{2}},0,0\right),\quad \vec{x}_4 =\left(\frac{1}{3},\frac{1}{\sqrt{3}},\frac{1}{3}\right).
    \]
    We can quickly exclude $\vec{x}_2$ by noting that $f'<0$ restricts $X<\frac{1}{2}$. To exclude $\vec{x}_3$ we note that if a trajectory converged to this point, then there exists a $T$ such that for all $t>T$, this trajectory must have $X<0$. But $X=a/h$, so $a<0$ eventually, forcing $\phi'<0$ eventually. This contradicts that $\phi\to\infty$ as $t\to\infty$.

    To exclude $\vec{x}_1$, we compute that the linearization of \eqref{eq:norm_system} at $\vec{x}_1$  has eigenvalues $-\frac{1}{2}$ and $\frac{1}{2}$ corresponding to a local stable manifold $Z=0$ and local unstable manifold $X=0$. This implies that no forward trajectory starting in $\mathcal{D}$ may converge to the point $\vec{x}_1$. 

    We now exclude periodic orbits in the interior of $\mathcal{D}$. The Dulac multiplier 
    \[
        B(X,Z) = \frac{1}{ZY^2} = \frac{2}{Z(1-2X^2-Z^2)}
    \]
    satisfies that
    \[
        \frac{\partial}{\partial X}(B\dot{X})+\frac{\partial}{\partial Z}(B\dot{Z}) = -\frac{2}{Z}<0
    \]
    so by the Bendixon-Dulac criterion the system has no closed orbit contained in the interior of $\mathcal{D}$.

    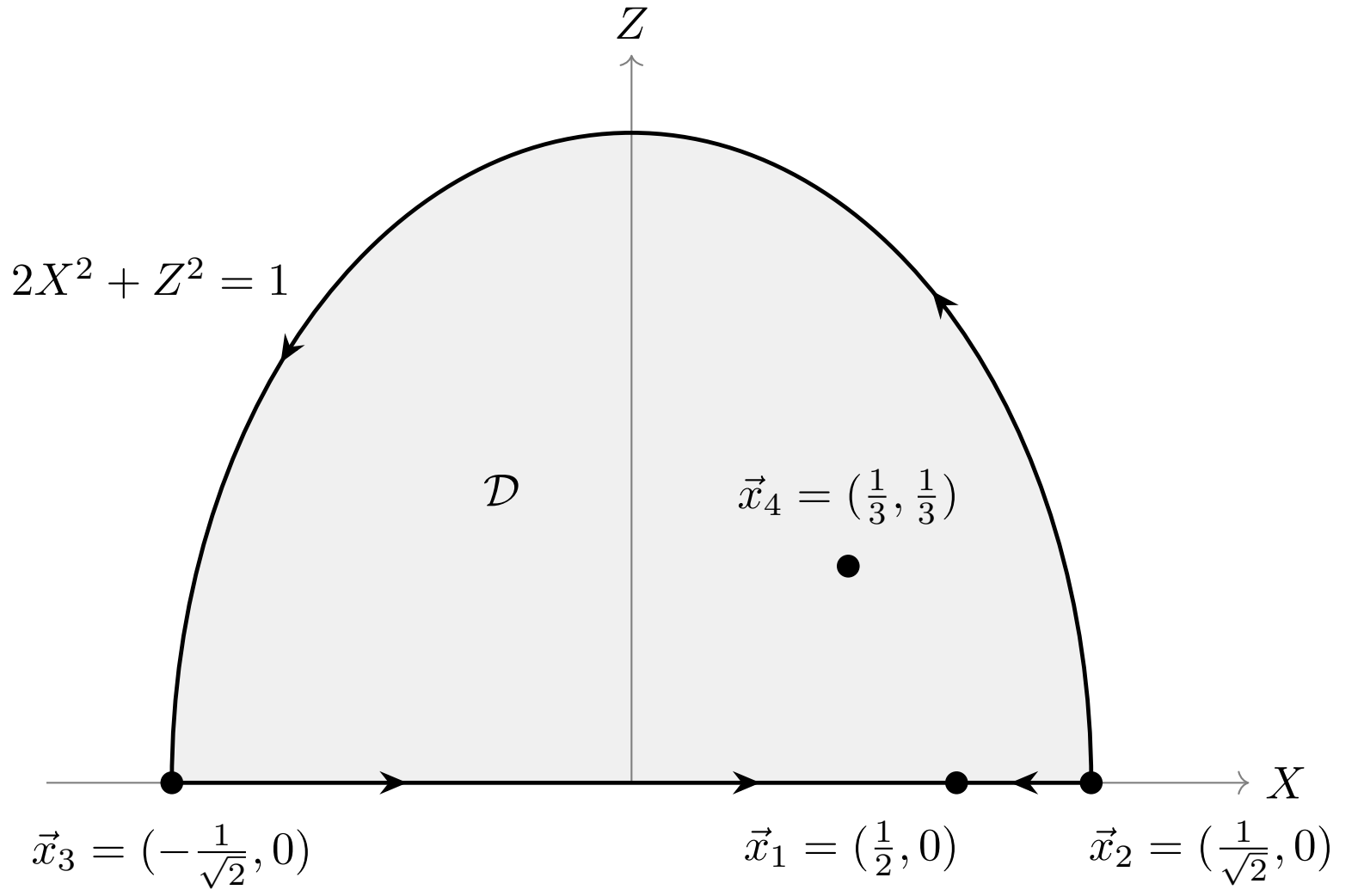
\begin{figure}
    \centering
    \begin{tikzpicture}[
        x=5cm, y=5cm,
        boundary/.style={line width=0.9pt},
        eq/.style={circle, fill, inner sep=0pt, minimum size=5pt},
        arrowat/.style={postaction={decorate, decoration={markings,
                        mark=at position #1 with {\arrow{Stealth[length=6pt, width=5pt]}}}}},
      ]
      \def\r{0.70711}   
    
      \fill[black!6] (\r,0) arc[x radius=\r, y radius=1, start angle=0, end angle=180] -- cycle;
    
      \draw[->, black!50] (-0.9,0) -- (0.95,0) node[right, black] {$X$};
      \draw[->, black!50] (0,0) -- (0,1.12) node[above, black] {$Z$};
    
      \draw[boundary, arrowat=0.3, arrowat=0.75]
        (\r,0) arc[x radius=\r, y radius=1, start angle=0, end angle=180];
      \draw[boundary, arrowat=0.3, arrowat=0.75] (-\r,0) -- (0.5,0);
      \draw[boundary, arrowat=0.6] (\r,0) -- (0.5,0);
    
      \node[eq] (P1) at (0.5,0) {};
      \node[eq] (P2) at (\r,0) {};
      \node[eq] (P3) at (-\r,0) {};
      \node[eq] (P4) at (1/3,1/3) {};
    
      \node[below left=3pt and -6pt of P1] {$\vec{x}_1=(\tfrac12,0)$};
      \node[below right=3pt and -6pt of P2] {$\vec{x}_2=(\tfrac1{\sqrt2},0)$};
      \node[below=3pt of P3] {$\vec{x}_3=(-\tfrac1{\sqrt2},0)$};
      \node[above=3pt of P4] {$\vec{x}_4=(\tfrac13,\tfrac13)$};
    
      \node[above left] at (-0.5,0.72) {$2X^2+Z^2=1$};
      \node at (-0.2,0.45) {$\mathcal D$};
    \end{tikzpicture}
    \caption{The region $\mathcal{D}$ with equilibria and the directions of the flow of the system along the invariant boundary.}
    \label{fig:D_with_flow}
    \end{figure}
    
    It remains to exclude possible $\omega$-limit points and polycycles in the boundary of $\mathcal{D}$. At $\vec{x}_2$, the linearization of \eqref{eq:norm_system} has eigenvalues $2-\sqrt{2}$ and $\sqrt{2}$, so it is a hyperbolic source (and thus not part of any polycycle). The point $\vec{x}_3$ has eigenvalues $2+\sqrt{2}$ and $-\sqrt{2}$, corresponding to the unstable manifold $Z=0$ and stable direction along the $Z$-direction respectively. Finally, at $\vec{x}_4$ the linearization has eigenvalues $-\frac{1}{3}\pm \frac{\sqrt{5}}{3}i$, hence it is a stable focus and cannot be part of a polycycle. Note that the boundary of $\mathcal{D}$ is invariant under the flow of \eqref{eq:norm_system} (since it corresponds to $Y=0$), so any polycycle must lie entirely in the boundary. We can rule out any boundary polycycle by inspecting the direction of the flow in Figure \ref{fig:D_with_flow}. Thus $\vec{x}_4$ is the only $\omega$-limit point of the system in $\mathcal{D}$, completing the proof.
\end{proof}

\begin{prop}\label{prop:crit_cone_asympt}
    When $\lambda=0$ we have the asymptotics as $t\to\infty$ satisfying
    \[
        h(t)=\frac{3}{t}+o(t^{-1}),\qquad \phi(t)=\frac{t}{\sqrt{3}}+o(t), \qquad e^{f(t)}=\frac{1}{t}+o(t^{-1}).
    \]
    In particular, we have
    \[
        g(t)\to dt^2+\frac{1}{3}t^2g_{S^2}
    \]
    so the soliton metric is asymptotic to a cone.
\end{prop}
\begin{proof}
    In Lemma \ref{lem:crit_asymptotics} we proved that as $t\to\infty$, 
    \[
        \frac{\phi'/\phi}{h}\to\frac{1}{3},\qquad \frac{1/\phi}{h}\to\frac{1}{\sqrt{3}}, \qquad \frac{e^f}{h}\to\frac{1}{3}.
    \]
    Since 
    \[
        \left(\frac{1}{h}\right)' = -\frac{h'}{h^2} = 2X^2+Z^2\to \frac{1}{3},
    \]
    we have that $th(t)\to 3$ as $t\to\infty$. Multiplying the three limits above by $t\,h(t)$ then yields
    \[
        t\frac{\phi'}{\phi}\to 1, \qquad \frac{t}{\phi}\to\sqrt{3}, \qquad te^f\to 1
    \]
    completing the proof.
\end{proof}

\subsection{Supercritical/compact solutions}
Finally, we show what happens in the case where $\lambda > 0$.  In this case, the solution to the system (\ref{eq:phi_f_sys}) does not exist for all time, and hence has maximal existence time at some $t = L$.  First, we prove this assertion.  Then, we show that for $\lambda \neq 2$, the sectional curvature diverges near this limiting point $t=L$, and there, no smooth compact closure can exist. If $\lambda = 2$, then we do in fact obtain a smooth closure, which is the Bismut flat $S^3$. 

We first show finite-time extinction of the system of ODEs (\ref{eq:phi_f_sys}).

\begin{lem}\label{lem:finite_ext_time}
    Let $\lambda > 0$, and let $(0, L)$ be the maximal interval on which the solution to the system of ODEs (\ref{eq:phi_f_sys}) exists and $\phi>0$.  Then
        \[ L \leq \pi \sqrt\frac{28}{3\lambda} < \infty. \]
    In particular, the solution cannot define a complete, noncompact, rotationally invariant gradient generalized Ricci soliton on $\mathbb{R}^3$.
\end{lem}

\begin{proof}
    Define the scalar quantity
        \[ W := 3h - 2a = 4\frac{\phi'}{\phi} - 3f'.\]
    We first compute the derivative
        \begin{align*}
            W' &= 3(-2a^2-c^2) - 2(b^2-c^2-ah)
            \\&= -6a^2-2b^2-c^2+2ah
             \\&= -\lambda - (4a^2-2ah+h^2).
        \end{align*}
    Next, we observe that
        \[ 4a^2-2ah+h^2 = \frac{3}{28}(3h - 2a)^2+\frac{1}{28}(10a-h)^2.  \]
    We do this to obtain the following bound on $W'$.  We now have 
        \begin{align} \label{eq:w_prime}
         W' = -\lambda - \frac{3}{28}W^2-\frac{1}{28}(10a-h)^2 \leq -\lambda - \frac{3}{28}W^2. 
         \end{align}
    Next, define
        \[ \theta(t) = \mathrm{arctan}\left (\sqrt \frac{3}{28\lambda} W(t)\right )  \]
    \noindent which is defined since $\lambda > 0$ by assumption.  We then compute
        \[ \theta' = \sqrt \frac{3}{28\lambda}\frac{W'}{1+\frac{3W^2}{28\lambda}} \leq -\sqrt \frac{3\lambda}{28}. \]
    Using the definition of $a, h$ and $W$, and using the Taylor expansions listed in \ref{eq:orig_asymp}, we can obtain the following expansions:
        \[
            a(t) = \frac{1}{t}+O(t), \qquad h(t) = \frac{2}{t} + O(t), \qquad W(t) = \frac{4}{t}+O(t).
        \]
    It clearly follows that
        \[ \displaystyle \lim_{t \to 0} \theta(t) = \frac{\pi}{2}.  \]
    Finally, integrating the expression for $\theta'$ from some $\epsilon > 0$ to some constant $t_0$ yields
        \begin{align*}
            \int_{\epsilon}^{t_0} \theta' \, dt \leq -\sqrt \frac{3\lambda}{28}\int_{\epsilon}^{t_0}dt.
        \end{align*}
    Evaluating the definite integral on both sides and taking the limit as $\epsilon \to 0$ yields
        \[ \theta(t_0) \leq \frac{\pi}{2} - \sqrt{\frac{3\lambda}{28}}t_0.  \]
    By definition of $\theta$, we also have 
        \[ -\frac{\pi}{2} \leq \theta(t_0). \]
    \noindent Putting these inequalities together and rearranging, we finally obtain
        \[ t_0 < \pi \sqrt {28/(3\lambda)}\]
    \noindent for any $t_0 \in (0, L)$.  This proves the result.
\end{proof}

\begin{thm}\label{thm:supercrit_nonsmooth_closures}
    Let $\lambda > 0$ and $(0,L)$ be the maximal interval on which $\phi > 0$ for the system (\ref{eq:phi_f_sys}). If $\lambda \neq 2$, then the metric $g=dt^2+\phi^2(t) g_{S^2}$, where $g_{S^2}$ is the standard round metric on $S^2$, has unbounded sectional curvature as $t$ approaches $L$. Therefore, no compact closure exists in this case. If $\lambda = 2$, then the solution is
        \[ \phi(t)=\sqrt{2}\sin(t/\sqrt 2), \qquad f(t) = 0 \]
    where adding the second pole gives the Bismut flat $S^3$ of radius $\sqrt{2}$.
\end{thm}

\begin{proof}
    First, we consider the case when $\lambda \neq 2$. Consider the function
        \[ J = \phi e^{-2f} \left ( \phi' - \frac{\phi h}{4}(2(\phi')^2+\phi^2e^{2f}) \right ).\]
    Differentiating, we obtain
        \[ J' =  e^{-2f} \left (1-(\phi')^2-\frac{1}{2}\phi^2e^{2f} \right )^2 \geq 0. \]
    Note that $J(0) = 0$, and, using the Taylor expansions for $f$ and $\phi$, we obtain the following expansion for $J$:
        \[ J(t) = \frac{(\lambda - 2)^2}{180}t^5+O(t^7). \]
    Hence, we can conclude that for $\lambda \neq 2$, we have that $J(t) > 0$ for $t \in (0, L)$, and $J$ is nondecreasing. If there was a second smooth pole, then one would require $\phi \to 0, \phi' \to -1$, $f$ has a finite limit, and $f' \to 0$.  However, this would imply that $J \to 0$, which contradicts the fact that $J(t) > 0$ and is nondecreasing. This rules out a smooth, compact solution. 

    To show that the curvature is unbounded for $\lambda \neq 2$, we consider the sectional curvatures 

    \[ K_{\mathrm{rad}} = -\frac{\phi''}{\phi} \qquad K_{\mathrm{tan}} = \frac{1-(\phi')^2}{\phi^2}. \]

    By contradiction, suppose $K_{\mathrm{rad}}$ and $K_{\mathrm{tan}}$ were both bounded near $L$. First, assume $\phi(L) > 0$.  Then we have that

    \[ (f')^2-5\frac{\phi'}{\phi}f'=  K_{\mathrm{rad}} + 3 K_{\mathrm{tan}}-\lambda. \]

    Since we also know $\phi'/\phi$ is bounded, then this expression bounds $f'$ and hence $f$.  The solution therefore stays in a compact subset, extending the solution past $L$, which contradicts the maximality of $L$.

    Next, assume $\phi(L) = 0$. Then boundedness of $K_{\mathrm{tan}}$ forces $\phi'(L) = -1$. The soliton equation, near $L$ gives
    \[ f' = \frac{\phi}{\phi'}(e^{2f}- K_{\mathrm{rad}}- K_{\mathrm{tan}}) \leq C\phi \]
    because $\phi' < 0$. Integrating bounds $f$ from above, and then we are able to obtain the bound $|f'| \leq C\phi$, and so $f$ has a finite limit. This again means that $J \to 0$ as $t \to L$, which contradicts positivity and monotonicity of $J$.  

    Therefore, for $\lambda > 0$ and $\lambda \neq 2$, we have
        \[ \displaystyle \lim_{t \to L}(|K_{\mathrm{rad}} | + |K_{\mathrm{tan}}|) = \infty \]
    establishing the curvature blowup. 

    When $\lambda = 2$, one easily verifies that $f=0$ and $\phi = \sqrt{2}\sin(t/\sqrt2)$ solves (\ref{eq:phi_f_sys}) and (\ref{eq:cons_law}).
\end{proof}

\begin{rmk}
    The quantity $J(t)$ has a geometric interpretation. Given orthonormal $X$, $Y$, the Bismut sectional curvatures satisfy
    \[
        K^{\pm}(X,Y)=K(X,Y)-\frac{1}{4}|H(X,Y,\cdot)|^2 = K(X,Y)-\frac{1}{2}e^{2f}
    \]
    where the last equality is special to dimension 3. Define $\mathcal{R}=\phi^2 \left(K_\text{tan}-\frac{1}{2}e^{2f}\right)$, where $\mathcal{R}$ is a scaling invariant measure of the Bismut sectional curvature of each orbit tangent plane. The definition of $J'$ can then be read as
    \[
        J' = e^{-2f}\mathcal{R}^2.
    \]
\end{rmk}

\begin{proof}[Proof of Theorem 1.1] 
    Part (a) of Theorem \ref{thm:main} follows from Proposition \ref{prop:glob_exist}. Part (b) follows from Proposition \ref{prop:crit_cone_asympt}. Parts (c) and (d) follow from Theorem \ref{thm:supercrit_nonsmooth_closures}.
\end{proof}

\subsection{Uniqueness}
The proof of uniqueness for the family of solitons constructed in this section amounts to proving that the $SO(3)$ action on $M$ must always admit a fixed point. Assuming no fixed point, the isotropy along a geodesic normal to a principal orbit implies the existence of a warping factor $\phi>0$ on all of $\mathbb{R}$, and then a Riccati comparison argument derived from the soliton equations excludes the existence of such a $\phi$. 

Before proving the fixed point result, we state and prove the Riccati comparison result.

\begin{prop} \label{prop:riccati_comp}
    Let $w \in C^1({\mathbb{R}})$ satisfy
    \[ w' \leq A - Bw^2  \] where $A \geq 0$, $B>0$ are constants. Then $|w| \leq \sqrt{A/B}$. In particular, setting $A=0$, if $w' \leq - Bw^2 $, then $w \equiv 0$ identically.
\end{prop}

\begin{proof}
    Fix $t_0 \in \mathbb{R}$, and first suppose $w(t_0) > \sqrt{A/B}$. Define $v(s) = w(t_0-s)$ for $s \geq 0$. Then $v' \geq Bv^2-A$. Since $v(0)>\sqrt{A/B}$, we have $v(s) \geq v(0)$ for all $s \geq 0$. Putting these inequalities together gives 
        \[ v' \geq v^2 \left ( B - \frac{A}{v^2} \right ) \geq  v^2 \left ( B - \frac{A}{v(0)^2} \right ) = \delta v^2  \]
    \noindent for some $\delta > 0$, which follows from the fact that $v(0) > \sqrt{A/B}$. Finally, from the inequality $v(s) \geq v(0)$, we have 
        \[ \frac{1}{v(s)} \leq \frac{1}{v(0)} - \delta s \]
    \noindent which would imply $v(s)<0$ for sufficiently large $s$, a contradiction. For $w(t_0) < -\sqrt{A/B}$, an analogous argument works, using $v(s) = -w(t_0+s)$. This proves the desired bound.
\end{proof}

\begin{lem}\label{lem:dim3_fixedpoint}
    Let $(M^3,g,f,H)$ be a complete gradient steady generalized Ricci soliton with $H \not\equiv 0$, and suppose there exists an effective $SO(3)$-action preserving $g$ and $H$ which has two dimensional principal orbits. Then the action has a fixed point.
\end{lem}

\begin{proof}
    By contradiction, suppose the action has no fixed point. Let $\gamma: \mathbb{R} \to M$ be a complete unit speed geodesic on $M$ which is normal to a principal orbit at $p$. This geodesic stays normal to principal orbits as geodesics stay normal to Killing fields. Let $K=(G_p)^0$ be the identity component of the stabilizer of $p$, and since $K$ acts trivially in the normal direction, we have $\mathfrak{k}\subseteq \mathfrak{g}_{\gamma(t)}$. Since the Lie algebra $\mathfrak{so}(3)$ has no two dimensional subalgebra and since the action does not have a fixed point by assumption, the isotropy of the action along $\gamma$ stays one dimensional. 
    
    The group action therefore allows us to locally write the metric $g$ as $g=dt^2+g_t$, where $g_t$ is a smooth family of $G$-invariant metrics on $G/K$. Since $G/K=SO(3)/SO(2)\cong S^2$, the invariance of $g_t$ under $G$ implies $g_t=\phi^2(t)g_{S^2}$, for a smooth, positive function $\phi$ on $M$. We can define $\phi > 0$ on all of $\mathbb{R}$ by defining $F: \mathbb{R} \times G/K \to M$ given by $F(t, gK) = g\gamma(t)$, which is a local diffeomorphism everywhere sine the isotopy stays one dimensional by assumption. The pullback metric $F^*g = dt^2+g_t$ where $\phi(t)$ is now given by $g_t = \phi^2(t)g_{S^2}$ now defines $\phi > 0$ on all of $\mathbb{R}$. Furthermore, $(f,\phi)$ solves the soliton system of equation \eqref{eq:phi_f_sys}.

    We now show that no smooth positive function $\phi$ can solve (\ref{eq:phi_f_sys}), obtaining the required contradiction. We split this up into three cases, based on the value of the conserved quantity $\lambda$.

    First, if $\lambda > 0$, then Lemma \ref{lem:finite_ext_time} already proves an upper bound to the existence of solutions to \eqref{eq:phi_f_sys}. If $\lambda = 0$, then as shown in the proof of Lemma \ref{lem:finite_ext_time}, we have that $W' \leq -\frac{3}{28}W^2$. Using the Ricatti comparison lemma stated and proved above, this implies that $W \equiv 0$ identically. Then by \eqref{eq:w_prime}, we have that
        \[ W' = -\lambda - \frac{3}{28}W^2-\frac{1}{28}(10a-h)^2 = -\frac{1}{28}(10a-h)^2.\]
    Hence, $10a-h=0$. However, we also have $W = 3h-2a = 0$ which implies $a=h=0$. Then by \eqref{eq:h_square}, we have that $2b^2+c^2=0$, a contradiction.
    Finally, if $\lambda < 0$, we again recall the inequality
        \[ W' \leq -\lambda - \frac{3}{28}W^2. \]
    Again, applying the previous lemma yields
        \begin{align} \label{eq:w_bound}
         |W| \leq  \sqrt{\frac{-28\lambda}{3}}.
        \end{align}
    Next, by \eqref{eq:h_square}, we have that $h^2 > -\lambda$. We can also assume $h > 0$ without loss of generality by swapping $t$ with $-t$ if necessary by reparameterizing $\gamma$. To justify this further, under the reparameterization $\tilde{\gamma}(t) := \gamma(-t)$, the sign of $h$ swaps, $b,c$ remain unchanged, the sign of $a$ swaps, and these functions satisfy the same equations with the same $\lambda$. Also, from \eqref{eq:h_square}, the sign of $\lambda$ does not change, so we are justified in assuming $h > \sqrt{-\lambda}$. Equation \eqref{eq:h_square} also gives us $|a| \leq h/\sqrt{2}$. Hence, $W=3h-2a \geq (3-\sqrt{2})h > 0$. Putting this together with (\ref{eq:w_bound}), this bounds $h$ on $\mathbb{R}$.  Looking again at (\ref{eq:h_square}), a $h$ bound also bounds $a, b, c$. Further, using the definitions of $a, b, c, h$ and the system \eqref{eq:phi_f_sys}, we can obtain the following expressions the the derivatives of $a, b, c, h$:

    \begin{align*}
        a' &= b^2-ah-c^2
        \\b'&= -ab
        \\c'&= (2a-h)c
        \\h' &= -2a^2-c^2 
    \end{align*}
    
    (\ref{eq:phi_f_sys})to check that the derivatives of $a, b, c, h$ are also bounded.  Boundedness of $h$ together with the fact that $h' \leq 0$ gives us that $h$ has a finite limit at $-\infty$. We will call this constant $\kappa$. Hence, for any fixed $t_0$, the integral
        \[ \int_{-\infty}^{t_0} (2a^2+c^2) \  dt < \infty \]
    converges. The function $2a^2+c^2$ is uniformly continuous since $a, c$ and their derivatives are bounded. Hence, in order for the integral above to converge, the function $2a^2+c^2$ must tend to zero, and so $a$ and $c$ also tend to 0 as $t \to -\infty$.  By \eqref{eq:h_square}, we also must have that
        \[ B^2 \to \frac{\kappa^2+\lambda}{2} :=D > 0. \]
    Also, we have that
        \[ a' = b^2-ah-c^2 \to D  \]
    as $t \to -\infty$, using the expression for $a'$ stated above. Finally, choose $T$ such that $a' > D/2$ for $t \leq T$. Then
        \[ a(T) - a(t) \geq \frac{D}{2}(T-t). \]
    The right hand side diverges as $t \to -\infty$, contradicting boundedness of $a$. 

    Hence, no such $\phi > 0$ can solve the system, contradicting the claim of no fixed point. Hence, the $SO(3)$ action must contain a fixed point as claimed.
\end{proof}

\begin{thm}\label{thm:3dim_uniqueness}
    Up to homothety, every complete, steady generalized Ricci soliton with $H \not\equiv 0$ admitting an effective, nontransitive $SO(3)$ action is a member of the classification of Theorem \ref{thm:main} or a quotient of the Bismut flat $S^3$.
\end{thm}
\begin{proof}
    First observe that the assumptions of the theorem imply that the principal orbits of the group action have dimension two. It was previously argued at the start of Lemma \ref{lem:dim3_fixedpoint} that if the action does not admit a fixed point, the isotropy of a point on a principal orbit must either be one- or zero-dimensional. If the isotropy is zero dimensional, then the $SO(3)$ action is transitive, which we are excluding here. Thus the isotropy must be one-dimensional, so the principal orbits are two-dimensional.
    
    Now by Lemma \ref{lem:dim3_fixedpoint}, any $SO(3)$-invariant steady generalized Ricci soliton with two-dimensional principal orbits admits a fixed point for the action. At the fixed point, the rotationally symmetric metric ansatz yields the soliton must satisfy the system of equations \eqref{eq:phi_f_sys}, and by rescaling we may fix the initial conditions for $\phi$ and $f$ as in the beginning of this section. Collecting the results of this section implies the soliton must be a member of the classification in Theorem \ref{thm:main}.
\end{proof}

\begin{rmk}
    To generalize the uniqueness proof to higher dimensions, where we consider an $SO(n)$ action for solutions of \eqref{eq:ngen_Ricci_sol}, we would need to require that the principal orbits are codimension one, as it is possible for the isotropy to jump up in dimension along a geodesic in this setting. Given that assumption, uniqueness can be proved via an essentially identical argument as in Theorem \ref{thm:3dim_uniqueness}.
\end{rmk}

\section{Higher Dimensional Generalized Bryant Solitons}
The generalized Ricci flow couples the Ricci flow to the evolution of a three-form, but one could also consider a general analogous coupling of the evolution of an $n$-form to the Ricci flow. This framework is investigated in even greater generality in \cite{streetsScalarCurvatureEntropy2023}, where $H$ may be an arbitrary closed section $\Lambda^* T^*M$. Furthermore, recent work motivates the coupling of $n$-forms to the Ricci flow from the string theory point of view \cite{papadopoulos2026braneflows}, and the methods of the previous section extend directly to this more general setting. The ODE analysis from the 3D generalized steady gradient Ricci solitons with $SO(3)$-symmetry translates essentially verbatim, with similar conclusions, to the $SO(n)$-symmetry ansatz for steady gradient $n$-generalized Ricci flow solitons. We state the necessary equations and sketch the results here.

Let $(M^n,g)$ be a Riemannian manifold with $H=db+H_0$ a closed $n$-form, for $n\geq 3$. We define the $n$-generalized Ricci flow by
\begin{align*}
    \frac{\partial}{\partial t} g &= -2(\Ric-c_n H^2),\\
    \frac{\partial}{\partial t} H &= -\Delta_{d,g_t} H
\end{align*}
where $H_{ij}^2=H_{ik_2\cdots k_n}H_j^{\;k_2\cdots k_n}$ and $c_n>0$ is a normalization constant. The steady gradient solitons for this flow satisfy
\begin{align}\label{eq:ngen_Ricci_sol}
    \Ric-c_n H^2+\nabla^2 f &= 0,\\
    d^*(e^{-f}H) &= 0.
\end{align}
Since $H$ is a top degree form, the second equation above yields $H=h_0 e^f dV_g$ for $h_0$ constant, so $H^2=(n-1)!h_0^2e^{2f}g$. Thus we have
\[
    \Ric+\nabla^2 f = \mu e^{2f}g
\]
where $\mu=c_n (n-1)!h_0^2$. The same arguments from the dimension three case allows us to fix $\mu=1$ and $f(0)=0$. We derive analogous soliton equations as in the generalized Ricci flow case by taking the trace and divergence of \eqref{eq:ngen_Ricci_sol} to yield the conservation law
\[
    \nabla(R+|\nabla f|^2-(2n-1)e^{2f})=0 \iff R+|\nabla f|^2-(2n-1)e^{2f}=\lambda.
\]

Assuming $SO(n)$-symmetry of the solitons, we have $g=dt^2+\phi(t)^2g_{S^{n-1}}$. Standard formulas (see \cite{petersonRiemannianGeometry2016}) yield sectional curvatures $K_{\text{rad}}=-\phi''/\phi$ and $K_{\text{tan}}=(1-\phi'^2)/\phi^2$, hence
\begin{align*}
    &\Ric(\partial_t,\partial_t)=(n-1)K_{\text{rad}},& &\Ric(e,e)=K_{rad}+(n-2)K_{\text{tan}},\\
    &\nabla^2 f(\partial_t,\partial_t)=f'',& &\nabla^2f(e,e)=\frac{\phi'}{\phi}f',
\end{align*}
where $e$ is a unit vector orthogonal to $\partial_t$. Substitution yields the $n$-dimensional system analogous to \eqref{eq:phi_f_sys}:
\begin{equation}\label{eq:n_phi_f_sys}
\begin{split}
    \phi'' &= (n-2)\frac{1-(\phi')^2}{\phi}+\phi'f'-\phi e^{2f},\\
    f'' &= (n-1)(n-2)\frac{1-(\phi')^2}{\phi^2}+(n-1)\frac{\phi'}{\phi}f'-(n-2)e^{2f}.
\end{split}
\end{equation}
The analogous initial conditions are $\phi(0)=f(0)=f'(0)=0$, $\phi'(0)=1$ and $\phi'''(0)=q$. The soliton conservation law then becomes
\begin{equation}\label{eq:n_cons_law}
    (n-1)(n-2)\frac{1-(\phi')^2}{\phi^2}+2(n-1)\frac{\phi'}{\phi}f'-(f')^2+e^{2f} = n(n-1)q+2n-1 \equiv \lambda,
\end{equation}
where $\lambda=0$ is the critical parameter and $\lambda_{\text{round}}=n-1$ yields the smooth soliton on $S^n$.

\subsection{Subcritical/noncompact solutions}

The proof of completeness when $\lambda\leq 0$ is essentially identical to the proof of Proposition \ref{prop:glob_exist}, so we only state the result here.

\begin{thm}[Subcritical completeness]
    If $\lambda\leq 0$, the local solution to \eqref{eq:n_phi_f_sys} exists on $[0,\infty)$ and yields a complete soliton on $\mathbb{R}^n$.
\end{thm}

The dilaton constant $\lambda$ determines the asymptotic curvature, growth rate of $f$ and radius of the profile function for the rotationally symmetric solitons. 

\begin{prop}[Complete, noncompact asymptotics]
    Let $\lambda<0$. As $t\to\infty$, we have that 
    \[
        h\to\sqrt{-\lambda},\qquad f'\to -\sqrt{-\lambda}, \qquad \phi^2 = \frac{2(n-2)}{\sqrt{-\lambda}}t+o(t).
    \]
    The soliton is asymptotically cylindrical in the sense that $tR\to \frac{(n-1)\sqrt{-\lambda}}{2}$.
\end{prop}
\begin{proof}
    The proof of the asymptotics for $h$ and $f'$ are the same as the proofs in Lemma \ref{lem:crit_h_phi_asymp}. Now consider the quantity $w=\phi\phi'$, and observe that
    \[
        w' = (n-2)-\phi^2 e^{2f}-w\left((n-3)\frac{\phi'}{\phi}-f'\right).
    \]
    Since $\frac{\phi'}{\phi}\to 0$ as $t\to\infty$, we have $(n-3)\frac{\phi'}{\phi}-f'\to\sqrt{-\lambda}$. This implies $w'=-(\sqrt{-\lambda}+o(1))w+(n-2)+o(1)$, yielding that $w\to \frac{n-2}{\sqrt{-\lambda}}$. This proves the asymptotics for $\phi^2$.

    To show the curvature asymptotics, we note $\phi'\to 0$, and the soliton equations \eqref{eq:n_phi_f_sys} imply $\phi''\to 0$. To leading order, this gives $R = \frac{(n-1)(n-2)}{\phi^2}$, completing the theorem.
\end{proof}

\subsection{Critical/conical solutions}
\begin{prop}
    At $\lambda=0$, the smooth solution is complete with asymptotics
    \[
        t\,h\to n,\quad \frac{\phi(t)}{t}\to\sqrt{\frac{n-2}{n}},\quad t\, e^{f}\to 1.
    \]
    In particular, the soliton is asymptotic to the cone $g=dt^2+\frac{n-2}{n}t^2g_{S^{n-1}}$.
\end{prop}
\begin{proof}
    The analogous results for the system \eqref{eq:n_phi_f_sys} follow from the same methods used to prove the results for \eqref{eq:phi_f_sys} in the critical case. One considers the normalized equation of \eqref{eq:n_cons_law} (cf. \eqref{eq:norm_cons_law}) , computes the evolution of each of the terms in the equation, eliminates the variable $Y$ to yield a planar system on a compact domain, and then applies the same ODE theory arguments to yield convergence to an interior equilibrium point.

    The interior equilibrium point is, in the substitution $X=\frac{a}{h}$, $Z=\frac{c}{h}$ as before,
    \[
        (X,Z) = \left(\frac{1}{n},\frac{1}{n}\right).
    \]
    and following the argument of Proposition \ref{prop:crit_cone_asympt} yields the theorem.
\end{proof}

\subsection{Supercritical/compact solutions}
The proof of finite extinction time in the $n$-dimensional case is identical to the proof in Lemma (\ref{lem:finite_ext_time}).
\begin{prop}[Finite extinction time]
    Let $\lambda>0$, and let $(0,L)$ be the maximal interval on which the solution to the system of ODEs \eqref{eq:n_phi_f_sys} exists and $\phi>0$. Then $L<\infty$.
\end{prop}

\begin{rmk}
    When $n=3$, Theorem \ref{thm:supercrit_nonsmooth_closures} proves that the only smooth compact soliton is the Bismut flat $S^3$. The proof method does not extend to the $n\geq 4$ case, and contrary to the $n=3$ case, numerical evidence suggests that there do exist smooth compact solitons for the associated $n$-form flow with $\lambda\ne (n-1)$. 
\end{rmk}

\printbibliography

\vspace{3mm}

\noindent(Eric Cochran) \textsc{Rowland Hall, Department of Mathematics, University of California, Irvine, Irvine, CA 92617} \\
\textit{E-mail: }\href{mailto:ecochra2@uci.edu}{ecochra2@uci.edu}

\vspace{3mm}

\noindent(Michael Womack) \textsc{Rowland Hall, Department of Mathematics, University of California, Irvine, Irvine, CA 92617} \\
\textit{E-mail: }\href{mtwomack@uci.edu}{mtwomack@uci.edu}

\end{document}

%% file: preamble.tex
\usepackage[utf8]{inputenc}
\usepackage{graphicx} 
\usepackage{amsfonts}
\usepackage{amsmath}
\usepackage{amsthm}
\usepackage{amssymb}
\usepackage{tikz}
\usetikzlibrary{arrows.meta,decorations.markings,positioning}
\usepackage{comment}
\usepackage{yhmath}
\usepackage{xcolor}
\usepackage[style=numeric,backend=biber]{biblatex}
\usepackage{hyperref}

\newtheorem{thm}{Theorem}
\newtheorem{prop}[thm]{Proposition}
\newtheorem{lem}[thm]{Lemma}

\newtheorem{rmk}[thm]{Remark}

\newtheorem{proposition}[thm]{Proposition}

\numberwithin{thm}{section}
\numberwithin{equation}{section}

\newcommand{\Ric}{\mathrm{Ric}}

\bibliography{bibfile}